\documentclass[12pt]{article} \textwidth=6in \oddsidemargin=0in
\usepackage{amssymb}
\begin{document}\def\ov{\over} \def\ep{\varepsilon}
\newcommand{\C}[1]{{\cal C}_{#1}} \def\inv{^{-1}}\def\be{\begin{equation}}
\def\ee{\end{equation}}\def\x{\xi}\def\({\left(} \def\){\right)} \def\iy{\infty} \def\e{\eta} \def\cd{\cdots} \def\ph{\varphi} \def\ps{\psi} \def\La{\Lambda} \def\s{\sigma} \def\ds{\displaystyle}
\newcommand{\xs}[1]{\x_{\s(#1)}} \newcommand\sg[1]{{\rm sgn}\,#1}
\newcommand{\xii}[1]{\x_{i_{#1}}} \def\ld{\ldots} \def\Z{\mathbb Z}
\def\a{\alpha}\def\g{\gamma}\def\b{\beta} \def\noi{\noindent}
\def\ar{\leftrightarrow} \def\d{\delta} \def\S{\mathbb S} \def\tl{\tilde} \def\E{\mathbb E} \def\P{\mathbb P} \def\TS{U\backslash T}\def\CS{{\cal S}} \newcommand{\br}[2]{\left[{#1\atop #2}\right]}
\def\bs{\backslash} \def\CD{{\cal D}} \def\N{\mathbb N} \def\t{\tau}
\def\la{\lambda} \def\R{\mathbb {R}} \def\ch{\raisebox{.4ex}{$\chi$}}
\def\sp{\vspace{1ex}} \def\G{\Gamma} \def\a{\alpha} \def\g{\gamma}
\def\tr{{\rm tr}\,} \def\dl{\delta}

\def\bc{\begin{center}} \def\ec{\end{center}}

\hfill  April 25, 2008
\begin{center}{\Large\bf A Fredholm Determinant Representation in ASEP}\end{center}
 
\begin{center}{\large\bf Craig A.~Tracy}\\
{\it Department of Mathematics \\
University of California\\
Davis, CA 95616, USA\\
email: \texttt{tracy@math.ucdavis.edu}}\end{center}

\begin{center}{\large \bf Harold Widom}\\
{\it Department of Mathematics\\
University of California\\
Santa Cruz, CA 95064, USA\\
email: \texttt{widom@ucsc.edu}}\end{center}
\begin{abstract}

In previous work \cite{TW} the authors found integral formulas for probabilities in the asymmetric simple exclusion process (ASEP) on the integer lattice $\Z$. The dynamics are uniquely determined once the initial state is specified. In this note we restrict our attention to the case of step initial condition with particles at the positive integers $\Z^+$ and consider the distribution function for the $m$th particle from the left. In \cite{TW} an infinite series of multiple integrals was derived for the distribution. In this note we show that the series can be summed to give a single integral whose integrand involves a Fredholm determinant. We use this determinant representation to derive (non-rigorously, at this writing) a scaling limit. 

\end{abstract}


\bc \textbf{I. Introduction}\ec

The asymmetric simple exclusion process (ASEP) is a basic  interacting particle model for nonequilibrium phenomena.
Since its introduction by Spitzer \cite{Spi} in 1970, it has become a popular and much studied
model. (See, e.g.\ \cite{Li1, Li2, PS, Se}.)  In ASEP on the integer lattice $\Z$ particles move according to two rules: (1) A particle at $x$ waits an exponential time with parameter one (independently of
all other particles), and then it chooses $y$ with probability $p(x,y)$; (2) If $y$ is vacant at that time it moves to $y$, while if $y$
is occupied it remains at $x$ and restarts its clock.  The adjective ``simple'' refers to the fact that allowed jumps are one
step to the right, $p(x,x+1)=p$, or one step to the left, $p(x,x-1)=1-p=q$.  Because we deal with continuous time, we need not worry about two or more particles attempting jumps at the same time.  The model is called T(totally)ASEP if either $p=1$ (particles hop only to
the right) or  $q=1$ (particles hop only to the left).

The dynamics are uniquely determined once we specify the initial state. We restrict our attention to the case of step initial condition with particles at the positive integers $\Z^+$.
With this initial condition it makes sense to talk about the position of the $m$th left-most particle at time $t$.  We denote this (random)
position by $x_m(t)$.  (So $x_m(0)=m$.)

In a now classic paper, Johansson \cite{Jo} for the case of TASEP with step initial condition showed that the distribution of
$x_m(t)$  is related to 
the distribution of the largest eigenvalue in the unitary Laguerre ensemble of random matrix theory.  This connection with
random matrix theory  leads to an expression  for
the distribution of $x_m(t)$
 as a Fredholm determinant.  The importance of the Fredholm determinant representation is that it makes possible an analysis of the regime  of  \textit{universal fluctuations}.  (See Johansson's Corollary 1.7.  For further discussion of this universal regime
see \cite{PS, Spo}.) 

For the  general ASEP model on $\Z$ there is,  as far as the authors know, no known connection
to random matrix theory and no analysis of the  analogous  scaling regime  as  analyzed by Johansson in the TASEP case.  
(However for stationary ASEP Bal\'azs and Sepp\"al\"ainen \cite{BS} and Quastel and Valk\'o \cite{QV} prove that the variance of the current across a characteristic is of order $t^{2/3}$ and that the diffusivity has order
$t^{1/3}$.)
In recent work \cite{TW} we showed (building on ideas from Bethe Ansatz \cite{Sc}) that $\P(x_m(t)\le x)$
can be expressed as an infinite series where the $k$th order term is a $k$-dimensional integral.  

In this paper we show that when $p$ and $q$ are nonzero this infinite
series can be summed to give a single integral whose integrand involves a Fredholm determinant. We use this determinant representation to derive a scaling limit for fixed $m$ with $x,\,t\to\iy$. (It is {\it not} an extension of the scaling limit of Johansson since $m$ is fixed.) It remains conjectural since the derivation lacks a justification for an interchange of limits, but it is surely true. 
\bc \textbf{II. Determinant Representation}\ec

The result is stated in terms of an operator with kernel
\[K(\x,\,\x')={\x^x\,e^{\ep(\x)t}\ov p+q\x\x'-\x},\]
where
\[\ep(\x)=p\,\x\inv+q\,\x-1.\]
If $q\ne 0$ and $\C{R}$ denotes the circle with center zero and radius $R$ then $K(\x,\,\x')$ is a smooth function on $\C{R}\times \C{R}$ when $R$ is sufficiently large. Then the  operator $K$ on $L^2(\C{R})$ defined by\footnote{All contour integrals are to be given a factor $1/2\pi i$.}
\[Kf(\x)=\int_{C_R}K(\x,\,\x')\,\,f(\x')\,d\x'\]
is trace class. With the notations
\[\t=p/q,\ \ \ \ (\la;\,\t)_m=(1-\la)\,(1-\la\,\t)\cdots(1-\la\,\t^{m-1}),\]
the result is that when $p$ and $q$ are both nonzero,
\be\P\left(x_m(t)\le x\right)=\int {\det(I-\la q K)\ov (\la;\,\tau)_m} {d\la\ov \la},
\label{Probrep}\ee
where the integral is taken over a circle so large that all the singularities of the integrand lie inside it. Evaluating the integral by residues gives the equivalent formula
\[\P\(x_m(t)>x\)=\sum_{i=0}^{m-1}\,\,{\det\(I-q\,\t^{-i}\,\,K\)\ov \prod\limits_{{j\ne i}\atop {j\le m-1}}\(1-\t^{j-i}\)}.\]
(When $\t=1$ this is modified in the obvious way.) In particular,
\[\P\(x_1(t)>x\)=\det(I-qK).\]

To state the Corollary to Theorem 5.2 of \cite{TW}, from which (\ref{Probrep}) will follow, we recall the the definition of the $\t$-binomial coefficient 
\[\br{N}{n}_\t={(1-\t^N)\,(1-\t^{N-1})\cdots (1-\t^{N-n+1})\ov (1-\t)\,(1-\t^2)\cdots (1-\t^n)}.\]
The result was that if the initial state is $\Z^+$ and $q\ne0$ then\footnote{The coefficients were given in \cite{TW} in terms of what we denoted there by $\br{N}{n}$, which is related to $\br{N}{n}_\t$ by $\br{N}{n}=q^{n(N-n)}\br{N}{n}_\t$.}
\[\P(x_m(t)=x)=(-1)^{m+1}\,\sum_{k\ge m}{1\ov k!}\,\br{k-1}{k-m}_\t\,p^{(k-m)(k-m+1)/2}\;q^{km+(k-m)(k+m-1)/2}\]
\[\times\int_{\C{R}}\cd\int_{\C{R}}\prod_{i\ne j}{\x_j-\x_i\ov p+q\x_i\x_j-\x_i}\;{1-\x_1\cd\x_k\ov\ds{\prod_i(1-\x_i)\,(q\x_i-p)}}\,\prod_i\(\x_i^{x-1}e^{\ep(\x_i)t}\)\,d\x_1\cd d\x_k,\]
where in the $k$-dimensional integral all indices run over $\{1,\ld,k\}$. 

If we sum this over $x$ from $-\iy$ to $x$ (which we may do since we may take $R>1$) we obtain
\[\P(x_m(t)\le x)=(-1)^{m}\,
\sum_{k\ge m}{1\ov k!}\,\br{k-1}{k-m}_\t\,p^{(k-m)(k-m+1)/2}\;q^{km+(k-m)(k+m-1)/2}\]
\be\times\int_{\C{R}}\cd\int_{\C{R}}\prod_{i\ne j}{\x_j-\x_i\ov p+q\x_i\x_j-\x_i}\;\prod_i{1\ov(1-\x_i)\,(q\x_i-p)}\,\prod_i\(\x_i^{x}e^{\ep(\x_i)t}\)\,d\x_1\cd d\x_k.\label{P}\ee

The observation that will allow us to express this in terms of Fredholm determinants is the identity
\be\det\left({1\ov p + q \x_i\x_j - \x_i}\right)_{1\le i,j\le k}=(-1)^k
 (pq)^{k(k-1)/2} \prod_{i\neq j}{\x_j-\x_i\ov p +q \x_i\x_j -\x_i}\,\,
\prod_i {1\ov (1-\x_i)(q\x_i-p)}.\label{xidet}\ee
This can be seen as follows. We may assume $\t\ne 1$ (i.e., $q\ne 1/2$) since both sides are continuous in $q$. If we make the substitutions
\[ \x_i={\e_i+1\ov \e_i +\t\inv}\, ,\]
then
\be{1\ov p+q \x_i \x_j -\x_i}=-{1\ov p(1-\t)}{(1+\t\e_i)\,(1+\t\e_j)\ov \e_i-\t \e_j}.\label{xieta}\ee
Since
\[\det\({1\ov\e_i-\t \e_j}\)\]
is a Cauchy determinant we can evaluate it, and we find that the determinant of (\ref{xieta}) equals
\be (-1)^k {\t^{k(k-1)/2}\ov p^k(1-\t)^{2k}}\,\prod_i{(1+\t\e_i)^2\ov \e_i}\,
\prod_{i \ne j}{\e_i-\e_j\ov \e_i-\t\e_j}.\label{etadet}\ee
We compute
\[{(1+\t\e_i)^2\ov \e_i}={p(1-\t)^2\ov (q\x_i-p)(1-\x_i)},\]
\[{\e_i-\e_j\ov \e_i-\t\e_j}=q\,{\x_j-\x_i\ov p+q\x_i\x_j-\x_i},\]
and find that (\ref{etadet}) equals the right side of (\ref{xidet}).

{}From (\ref{xidet}) we see that when also $p\ne0$ (\ref{P}) may be written
\pagebreak
\[\P(x_m(t)\le x)=(-1)^{m}\,
\sum_{k\ge m}{(-1)^k\ov k!}\,\br{k-1}{k-m}_\t\,p^{m(m-1)/2-k(m-1)}\;q^{-m(m-1)/2+km}\]
\[\times\int_{\C{R}}\cd\int_{\C{R}}\det(K(\x_i,\,\x_j))_{1\le i,j\le k}\;d\x_1\cd d\x_k\]
\[=(-1)^{m}\,\t^{m(m-1)/2}
\sum_{k\ge m}{(-1)^k\ov k!}\,\br{k-1}{k-m}_\t\,\({p\ov\t^m}\)^k\]
\be\times\int_{\C{R}}\cd\int_{\C{R}}\det(K(\x_i,\,\x_j))_{1\le i,j\le k}\;d\x_1\cd d\x_k.\label{intsum}\ee

The last integral is a coefficient in the Fredholm expansion of $\det(I-\la K)$. In fact
\be\sum_{k=0}^\iy{(-\la)^k\ov k!}\int_{\C{R}}\cd\int_{\C{R}}\det(K(\x_i,\,\x_j))_{1\le i,j\le k}\;d\x_1\cd d\x_k=\det(I-\la K).\label{exp}\ee
{}From \cite[p.26]{M} we have for $|z|$ small enough 
\[\sum_{k\ge m} \left[{k-1\atop k-m}\right]_\t z^{k}=z^m\sum_{j\ge 0} \left[{m+j-1\atop j}\right]_\t z^j=z^m\,\prod_{i=0}^{m-1}{1\ov  1- \t^i z} =\prod_{j=1}^{m}{z\ov  1- \t^{m-j} z}.\]
If we set $z=p\t^{-m}\,\la\inv$ this gives for $|\la|$ large enough
\be (-1)^m\,\t^{m(m-1)/2}\,\sum_{k\ge m} \left[{k-1\atop k-m}\right]_\t\({p\ov\t^m}\)^k\,\la^{-k}=\prod_{j=1}^{m}{1\ov 1-\la p\inv\t^{j}}={1\ov (\la/q;\,\t)_m}.\label{prod}\ee
If we multiply (\ref{exp}) and (\ref{prod}), multiply by $\la\inv$, and integrate and we see that (\ref{intsum}) is the same as (\ref{Probrep}).

\sp
\noi{\bf Remark}. Formula (\ref{Probrep}) holds when $p$ and $q$ are nonzero. When $p\to0$ we must get the TASEP determinant for the probability, and the question arises whether this is easy to see. The answer seems to be that it is not easy, but it can be derived from it with some work.
 
\bc \textbf{III. A Scaling Conjecture}\ec

Denote by $K_0$ the operator on $L^2(\R)$ with kernel
\[K_0(z,\,z')={1\ov\sqrt{2\pi}}\,\,e^{-(p^2+q^2)\,(z^2+{z'}^2)/4+pq\,zz'}.\]
(This is the symmetrization of the Mehler kernel.) When $p\ne q$ this is trace class. The conjecture is that if $p<q$, so there is a drift to the left, we have for each $m$
\sp
\be\lim_{t\to\iy}\P\(x_m(t)\le (p-q)t+(q-p)\,y\,t^{1/2}\)=\int {\det(I-\la q K_0\,\ch_{(-y,\iy)})\ov (\la;\,\tau)_m} {d\la\ov \la},\label{conjecture}\ee
\sp
where, as before, the integral is taken over a circle so large that all the singularities of the integrand lie inside it.

We know from (\ref{Probrep}) that
\[\P\(x_m(t)\le (p-q)t+(q-p)\,y\,t^{1/2}\)=\int {\det(I-\la q K)\ov (\la;\,\tau)_m} {d\la\ov \la},\]
where $K$ is as before and
\be x=(p-q) t+(q-p)\,y\,t^{1/2}.\label{x}\ee
Therefore the conjecture would follow if we can show that
\be \det(I-\la K)\to\det(I-\la K_0\,\ch_{(-y,\iy)})\ \textrm{as $t\to\iy$, uniformly on compact $\la$-sets}.\label{detlim}\ee

The Fredholm determinants are entire functions of $\la$, and the coefficients in their expansions about $\la=0$ are universal polynomials in the traces of powers of the operators. By considering the coefficients successively we can see that if (\ref{detlim}) is true then necessarily
\sp

\noi(i)\ {\it $\tr K^n\to\tr (K_0\,\ch_{(-y,\iy)}))^n$ as $t\to \iy$ for $n=1,\,2,\,\cd$\,.}
\sp 

If (i) holds then each coefficient in the expansion of $\det(I-\la K)$ converges to the corresponding coefficient in the expansion of $\det(I-\la K_0\,\ch_{(-y,\iy)})$. This is not sufficient to give (\ref{detlim}). It would be sufficient if we also knew that
\sp

\noi(ii)\ $\det(I-\la K)$\ {\it is uniformly bounded for large $t$ on compact $\la$-sets}.
\sp 

\noi For if this holds one sees from the Cauchy inequalities that the convergence of the series for $\det(I-\la K)$ is uniform in $t$.

We shall show that (i) holds.   

First, though, we remark that if the conjecture is true then we expect the right side of (\ref{conjecture}) to be a distribution function in $y$ for each $m$. As $y\to-\iy$ the numerator in the integrand approaches one, and expanding the contour shows that the resulting integral equals zero. But what about the $y\to+\iy$ limit, which should equal one? The right side becomes
\be 1-\sum_{i=0}^{m-1}\,\,{\det\(I-q\,\t^{-i}\,\,K_0\)\ov \prod\limits_{{j\ne i}\atop {j\le m-1}}\(1-\t^{j-i}\)},\label{detexp}\ee
where $K_0$ acts on all of $\R$. If this is to equal one when $m=1$ then we must have $\det(I-qK_0)=0$. If this holds and (\ref{detexp}) is equal to one when $m=2$ then we must also have $\det(I-q\t\inv K_0)=0$. And so on. The conclusion is that $\det\(I-q\,\t^{-i}\,K_0\)$ should equal zero for all $i\ge 0$, in other words that all $\t^i/q$ should be eigenvalues of $K_0$. These are indeed eigenvalues and the corresponding eigenfunctions are
\[e^{-(q^2-p^2)\,z^2/4}\,H_i\(\sqrt{{q^2-p^2\ov2}}\, z\),\]
where the $H_i$ are the Hermite polynomials.

Next, we explain where the conjecture came from: there is a kernel with the same Fredholm determinant\footnote{We use the term ``Fredholm determinant'' to mean the infinite series. In case the kernel is continuous and trace class on some Hilbert space this is the same as the operator determinant.} as $K$ that converges pointwise to a kernel with the same Fredholm determinant as $K_0\,\ch_{(-y,\iy)}$.
  
When $p<q$ we have, when $|\x|$ and $|\x'|$ are large enough,  
\[K(\x,\,\x')=(q-p){\x^x\,e^{\ep(\x)t}\ov(q\x-p)\,(q\x'-p)}\,\int_0^\iy e^{z\,\(q{1-\x'\ov q\x'-p}-p{1-\x\ov q\x-p}\)}\,dz.\]
Therefore if we define 
\[A(\x,\,z)=(q-p){\x^x\,e^{\ep(\x)t}\ov q\x-p}e^{-pz{1-\x\ov q\x-p}},\ \ \
B(z,\,\x')={1\ov q\x'-p}e^{qz{1-\x'\ov q\x'-p}},\]
then
\[K(\x,\,\x')=\int_0^\iy A(\x,\,z)\,B(z,\,\x')\,dz.\]
The kernel
\be K_1(z,\,z')=\int_{\C{R}}B(z,\,\x)\,A(\x,\,z')\,d\x=
{q-p\ov2\pi i}\int_{C_R}{\x^x\,e^{\ep(\x)t}\ov(q\x-p)^2}\,
e^{(qz-pz'){1-\x\ov q\x-p}}\,d\x\label{K1}\ee
on $\R^+$ has the same Fredholm determinant as $K$ on $\C{R}$.

The part of the exponent with the factor $t$,
\be (p-q)\,\log\x+\ep(\x),\label{tpart}\ee
has a critical point at $\x=1$, its second derivative is positive there, and on the line ${\rm Re}\,\x=1$ its real part has an absolute maximum there. When $x<0$ we may replace $\C{R}$ by this line, and we
make the substitution $\x\to 1+i\x t^{-1/2}$. We also make the
substitutions 
\[z\to (q-p)z t^{1/2},\ \ \ z'\to (q-p)z' t^{1/2},\]
and  multiply the kernel by $(q-p)\,t^{1/2}$, so the Fredholm determinant is unchanged. We let $t\to\iy$ and obtain the pointwise limit 
\[K_2(z,\,z')={1\ov\sqrt{2\pi}}e^{-(qz-pz'-(q-p)y)^2/2}.\]
Here $z,\,z'\in \R^+$. We may replace this by the kernel 
$K_3(z,\,z')\,\ch_{(-y,\iy)}(z')$,
where $z,\,z'\in\R$ and
\[K_3(z,\,z')={1\ov\sqrt{2\pi}}e^{-(qz-pz')^2/2}.\]

The kernel $K_0(z,\,z')$ is the symmetrization of $K_3(z,\,z')$. In fact
\[K_0(z,\,z')=e^{(q-p)z^2/4}\,K_3(z,\,z')\,e^{-(q-p){z'}^2/4},\]
so they have the same Fredholm determinants. Hence, the conjecture. 

Finally, we establish (i). Instead of $K$ we may use the kernel $K_1$ given by (\ref{K1}) because the traces of the powers are the same. We make the variable change  
\[\e={\x-1\ov q\x-p}\]
and find that with $x$ given by (\ref{x})
\[K_1(z,\,z')={1\ov2\pi i}\int_\g e^{-(qz-pz')\,\e+f(\e)t^{1/2}+ g(\e)\,t}\,d\e,\]
where
\[f(\e)=(q-p)y\,\log\({1-p\e\ov1-q\e}\),\]
\[g(\e)=(p-q)\,\log\({1-p\e\ov1-q\e}\)+(q-p)^2\,{\e\ov (1-p\e)\,(1-q\e)},\]
and $\g$ is a little circle around $1/q$ described clockwise. (The apparent complication here will make things  simpler later.) 

The function $g(\e)$ equals (\ref{tpart}) after our substitution. The line ${\rm Re}\,\x=1$ corresponds to the circle $\G$ with diameter $(0,\,1/q)$. On this contour $g$ has a critical
point at $\e=0$, its second derivative is positive there, and the real part of $g(\e)$ has an absolute maximum there.

We consider first 
\be\tr K_1={1\ov q-p}\,{1\ov2\pi i}\int_\g e^{f(\e)\,t^{1/2}+g(\e)\,t}\,{d\e\ov\e}.\label{trK1}\ee
In the neighborhood of $\e=0$
\[f(\e)=(q-p)^2y\e +O(\e^2),\ \ \ g(\e)=(q-p)^2\,\e^2/2+O(\e^3).\]
Since the singularity at $1/p>1/q$ is outside $\g$ and $\G$ we may deform $\g$ to $\G$ if we replace the denominator $\e$ in (\ref{trK1}) by $\e+0$.
We then make the substitution $\e\to  i\e t^{-1/2}$ and find (recalling that $\g$, and so $\G$, is described clockwise) that
\[\lim_{t\to\iy}\tr K_1={1\ov q-p}\,{1\ov2\pi}\int_{-\iy}^\iy e^{-(q-p)^2\e^2/2+i(q-p)^2y\e}\,{d\e\ov i\e+0}\]
\[={q-p\ov2\pi}\int_{-\iy}^\iy \int_0^\iy e^{-(q-p)^2\e^2/2+i(q-p)^2y\e}\,e^{-i(q-p)^2z\e}\,dz\,d\e
={1\ov\sqrt{2\pi}}\int_0^\iy e^{-(q-p)^2(z-y)^2/2}\,dz\]
\[={1\ov\sqrt{2\pi}}\int_{-y}^\iy e^{-(q-p)^2z^2/2}\,dz
=\tr K_3\,\ch_{(-y,\iy)}=\tr K_0\,\ch_{(-y,\iy)}.\]

Now comes the tricky part. We compute that
\[\tr K_1^n={1\ov (2\pi i)^n}\int_\g\cd\int_\g {e^{\,\sum_j (t^{1/2}\,f(\e_j)+t\,g(\e_j))}\ov \prod_j(q\,\e_j-p\,\e_{j+1})}\,d\e_1\cd d\e_n,\]
where we set $\e_{n+1}=\e_1$. Recalling that $\t=p/q$ we may write this as
\be{1\ov (2\pi iq)^n}\int_\g\cd\int_\g {e^{\,\sum_j (t^{1/2}\,f(\e_j)+t\,g(\e_j))}\ov \prod_j(\e_j-\t\,\e_{j+1})}\,d\e_1\cd d\e_n.\label{trint}\ee

We want to deform all contours $\g$ to $\G$. Suppose we deform the $\e_j$-contour $\g$, while perhaps some are still $\g$ and others are already $\G$. Since $\t<1$ the pole at $\e_j=\t\,\e_{j+1}$ will be crossed in the deformation whether $\e_{j+1}$ is on $\g$ or $\G$, while the pole at $\t\inv \e_{j-1}$ will not be crossed. The residue at the crossed pole will give rise to a lower-order integral. In the deformation of all contours to $\G$ we will eventually get a sum of integrals over $\G$ of order $\le n$. But there is the problem that the integrands in these integrals will be singular at zero, so we cannot blithely deform all contours to $\G$, even if we keep track of the residues.

What we do is first expand $\G$ slightly to a circle $\G'$ with diameter $(-\a,\,1/q)$, where $\a$ is very small and positive. Because $\G$ is expanded, but just a little, we will still cross poles at the $\t\e_j$ (for $\e_j$ on the  contour $\g$ or the entire contour $\G'$) but not at the $\t\inv \e_j$. (If we had shrunk $\G$ instead we would cross the pole at $\t\e_j$ only for some $\e_j\in\G'$, which would complicate determining the lower-order integrals.)

Expanding all contours as described we get a sum of integrals over $\G'$ of order $\le n$, one for each subset $S$ of $\{1,\ld,n\}$, in which the integrands are nonsingular. The integral corresponding to $S$ is obtained as follows: for each $j\in S$ remove the factor $\e_j-\t\e_{j+1}$ from the denominator, make the replacement $\e_j\to\t\e_{j+1}$ in the remainder of the integrand, and multiply by $2\pi i$. The integral is taken with respect to the $\e_j$ with $j\not\in S$. This is the integrand we get after we take the residues that arise from crossing the poles at the $\e_j=\t\e_{j+1}$ with $j\in S$.\footnote{Although $\g$ is expanded outward we multiply by $2\pi i$ rather than $-2\pi i$ because $\g$ is described clockwise.}

The result may be described equivelently as follows: replace the product
\[\prod_j(\e_j-\t \,\e_{j+1})\inv\] 
with all $\e_j\in \g$ by
\[\prod_j\Big[(\e_j-\t \,\e_{j+1})\inv+2\pi i\dl(\e_j-\t \,\e_{j+1})\Big]\]
with all $\e_j\in \G'$.
It follows that (\ref{trint}) is equal to
\[{1\ov (2\pi iq)^n}\int_{\G'}\cd\int_{\G'}e^{\,\sum_j (t^{1/2}\,f(\e_j)+t\,g(\e_j))}\,\prod_j\Big[(\e_j-\t \,\e_{j+1})\inv+2\pi i\dl(\e_j-\t \,\e_{j+1})\Big]\,d\e_1\cd d\e_n.\]

Once we have this we can take the $t\to\iy$ limit. Recall that in the above $\a$ could have been any sufficiently small positive number. We take it to be $t^{-1/2}$, so the left-most point of $\G'$ is $-t^{-1/2}$, and $\G'$ is vertical there. On $\G'$ we make the substitutions $\e_j\to i\e_j\,t^{-1/2}$ and obtain in the $t\to\iy$ limit
\[{1\ov (2\pi iq)^n}\int_{-\iy+i}^{\iy+i}\cd\int_{-\iy+i}^{\iy+i}\prod_j e^{-(q-p)^2\e_j^2/2+i(q-p)^2y\e_j}\]
\be\times\prod_j\Big[(\e_j-\t\e_{j+1})\inv+2\pi i \dl(\e_j-\t\e_{j+1})\Big]\,d\e_1\cd d\e_n\label{int1}\ee
We used here the fact that $\dl$ is homogeneous of degree $-1$.

We now undo what we did before. Suppose we had the integral
\be{1\ov (2\pi iq)^n}\int_{-\iy-i}^{\iy-i}\cd\int_{-\iy-i}^{\iy-i}\prod_j e^{-(q-p)^2\e_j^2/2+i(q-p)^2y\e_j}\,\prod_j(\e_j-\t\e_{j+1})\inv\,d\e_1\cd d\e_n,\label{int2}\ee
and wanted to integrate over Im$\,\e_j=+1$ instead over Im$\,\e_j=-1$. If we raised the $\e_j$-contours successively we would obtain a sum of integrals of lower order over Im$\,\e_j=1$ that come from the residues at various poles at $\e_j=\t\e_{j+1}$. The sum of all integrals that arise is precisely equal to (\ref{int1}). Thus the limit of $\tr K^n$ is equal to (\ref{int2}).

We now substitute into (\ref{int2}) the integral representations, valid when Im$\,\e_j=-1$,
\[{1\ov i(\e_j-\t\e_{j+1})}=\int_0^\iy e^{-i(\e_j-\t\e_{j+1})z_j}\,dz_j,\]
and integrate first with with respect to the $\e_j$. We obtain 
\[{1\ov(2\pi)^{n/2}}\int_0^\iy\cd\int_0^\iy\prod_{j=1}^n e^{-(qz_{j-1}-pz_j-(q-p)y)^2/2}\,dz_1\cd dz_n,\]
where $z_{0}=z_n$. This is equal to
\[{1\ov(2\pi)^{n/2}}\int_{-y}^\iy\cd\int_{-y}^\iy\prod_{j=1}^n e^{-(qz_{j-1}-pz_j)^2/2}\,dz_1\cd dz_n=\tr (K_3\,\ch_{(-y,\iy)})^n=\tr (K_0\,\ch_{(-y,\iy)})^n.\]

Thus (i) is established. What would it take to establish (ii)? It would be enough to find kernels that have the same Fredholm determinants as $K$ and that have bounded trace norms on some Hilbert space. We could use the Hilbert-Schmidt norm instead (which is weaker and easier to compute) since that would give uniform boundedness of the regularized 2-determinant, and since the traces are bounded the determinant would also be. But so far we have not found any such kernels.

\end{document}